# Application of hybrid classical-quantum annealing technology to the 3D Bin-Packing Problem

Mohsen Rahmani[a]     Nitish Umang[b]     Mason Christensen[a]

[a] D-Wave Quantum Inc., Professional Services, Burnaby, British Columbia, Canada

[b] Johnson & Johnson Supply Chain Digital & Data Science, USA

## Abstract

In this paper, we study the use of hybrid classical-quantum annealing technology to solve a critical business optimization problem that is a form of multiproduct multi-bin packing in three dimensions with support constraints and case orientations along all three axes. We developed an exact mathematical model based on mixed-integer programming (MIP) to solve the problem, using fewer variables than previously existing models. Furthermore, to ensure the stability of the cases within bins, the model employs a novel formulation to represent the support constraints. We then compared and analyzed the solution performance of the classical solver Gurobi and D-Wave's constrained quadratic model (CQM) solver on the MIP model, both with and without support constraints. Results from the computational studies offer valuable insights into how the hybrid classical-quantum solver compares against widely used purely classical solvers on computationally hard problems, such as the problem of interest. This comparison examines aspects such as the optimality gap and volume utilization across different computational time limits.

# 1. Introduction

Quantum computing (QC) is a rapidly growing technology that uses concepts from quantum mechanics to solve many real-world problems using specialized hardware. Due to the current limitations in QPU hardware, researchers employ hybrid classical-quantum algorithms to address large-scale problems, where portions of the problem are algorithmically reduced to sizes that the QPU can manage. In recent years, annealing-based quantum computing has found applications in several domains, including machine learning (Biamonte et al., 2017; Khoshaman et al., 2018; Amin et al., 2018; Gircha et al., 2023; Ajagekar & You, 2023), financial modeling (Mugel et al., 2022), materials and molecular design (Gircha et al., 2023; Ajagekar & You, 2023), scheduling (Ikeda et al., 2019; Venturelli et al., 2015) , vehicle routing (Feld et al., 2019), and telecommunications (Kim et al., 2019), among others.  We refer the reader to Yarkoni et al. (2022) for a broad literature review on the industry applications of quantum annealing-based optimization.

Quantum annealing (QA) is a quantum computing approach especially well-suited for solving optimization problems (Kadowaki et al.,1998), including satisfiability and Ising model problems. QA belongs to the adiabatic quantum computing paradigm of computation (Farhi et al., 2000), which guarantees that a quantum system tracks the low-energy state(s) of an evolving Hamiltonian H(t) as long as the time-evolution is sufficiently slow (D-Wave, 2024). This principle is applied in D-Wave's quantum annealing  processors, where a time-dependent Hamiltonian is used to evolve the system from an initial equal-superposition state to a final state encoding low-energy solutions of the problem; King et al. (2024) demonstrated computational supremacy in quantum simulation in superconducting quantum annealing, and  both Lanting et al. (2014) and

King et al. (2022) have demonstrated quantum coherence and entanglement in QA processors, suggesting their viability for large-scale quantum computing. QA has shown promise in optimization and simulation applications, with potential advantages over classical methods in terms of solution quality and computation time for certain problems when used directly as a solver, or within a hybrid framework (McGeoch & Farré, 2023; Raymond et al., 2023).

The primary focus of this paper is to propose a new mathematical formulation for the 3D bin-packing problem that could be solved using a hybrid solver solution provided by D-Wave. The secondary objective is to evaluate the performance of this solver in comparison to a well-established solver recognized by the operations research community, in terms of solution quality, time to solution and tractability. In this paper, we do not explore the specifics of the solvers, the contribution of the QPU within the solution method, the problem embedding process, or the role of the quantum computer in the hybrid methods, as these details have not yet been disclosed.

The 3D bin-packing problem has been extensively studied due to its practical importance and computational challenges. Two main approaches exist: heuristic methods and exact methods. Heuristic methods, such as first fit decreasing (FFD) and best fit decreasing (BFD) (Martello & Toth, 1990) sequentially add items to bins. Advanced heuristics, like extreme point-based heuristics (Crainic et al., 2008), and maximal space strategies (Parreño et al., 2008), as well as global search approaches, such as biased random key genetic algorithm (Gonçalves & Resende, 2013) and guided local search (Faroe et al., 2003), combine multiple strategies to improve solution quality. While heuristics are effective for large-scale problems, exact methods provide precise solutions and benchmarks for smaller instances. Martello et al. (2000) presented exact algorithms for single and multi-bin 3D bin-packing problem, which have since been the foundation for many subsequent studies. Junqueira et al. (2012) added practical constraints, such as cargo stability and load-bearing limitations to this problem. Paquay et al. (2016) proposed the MIP model for application to air cargo loading, including a wide array of practical constraints, such as cargo stability and load-bearing constraints.

This study explores the application of hybrid classical-quantum annealing computing to solve the 3D bin-packing problem. We propose a new formulation that reduces the number of variables needed to model geometric constraints and the placement of the cases within multi-bins. Additionally, we propose a new formulation that models support constraints to ensure stability of the cargo within a bin. Finally, we present results demonstrating both the potential advantages and current limitations of applying hybrid quantum annealing to the 3D bin-packing problem. The remainder of the paper is organized as follows: in section 2, we present a mathematical model for the 3D bin-packing problem, introducing our novel approach to modeling this problem for multiple bins, case orientations, and support constraints. Detailed computational comparison between the quantum-based constrained quadratic model (CQM) solver and the classical Gurobi solver is presented in section 3. In section 4, we conclude the paper with key findings from the computational study and our outlook with respect to the use of quantum hybrid algorithms to solve optimization problems. The appendices include detailed descriptions of the mathematical models, complete information on the test cases, and further analysis of the results.

## 2. Mathematical Models for 3D Bin-Packing Problem

In this section, we describe the problem and the solution methods for our specific problem of interest. The problem is a form of multiproduct 3D bin packing with side operational constraints and presents itself in several different forms across different business use cases. It may involve loading cases into shipper boxes, gaylords, or pallets, stacking gaylords on trucks, loading pallets into ocean shipping containers, and so on.

One of the key challenges is obtaining near-optimal or high-quality solutions as the problem scales, i.e., as the count and heterogeneity in cases and bins increases. Classical approaches based on exact methods like mixed integer programming typically yield global optimal or near-optimal solutions for small problem sizes but are either too slow or break down and do not return any results for realistic problem sizes. Metaheuristics, conversely, are computationally much faster, but the solution quality is often less than satisfactory. From a technical standpoint, the goal of this study is to examine the potential benefits and shortcomings of hybrid quantum methods over classical optimization in terms of time to solution, solution quality, and problem scaling.

### 2.1 Problem Definition

The 3D bin-packing problem, as defined for the purpose of this study, borrows elements from several business use cases related to bin packing. The primary objective of the resulting optimization problem is to load cases into cuboidal bins to maximize space utilization. Improved space utilization leads to reduced number of bins, which in turn reduces transportation costs

(including fuel cost and CO2 emissions) and relieves congestion, thereby alleviating pressures on the supply chain. Other ancillary objectives may include:

- Maximizing homogeneity in bin packing with respect to product distribution, i.e., as far as possible, cases belonging to the same product category must be loaded together in the same bin(s).
- Optimizing the distribution of bin sizes, where ideally a larger bin (say, a 32 L shipper box container) is preferred over two or more smaller bin (say, two 16 L shipper boxes) for the same volume to be loaded.

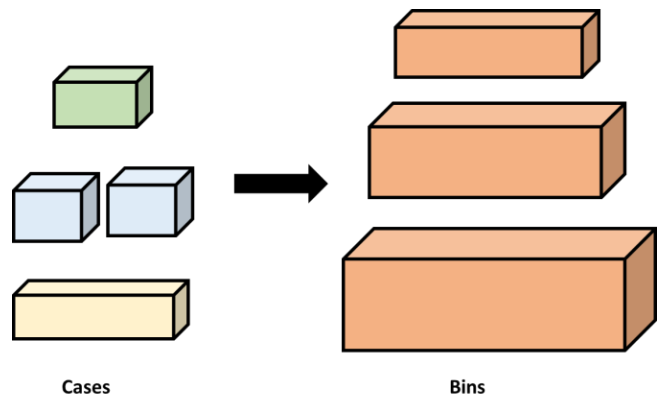

Figure 1. Loading multiproduct cases into bins of different sizes.

In the base model, the arrangement of cases in each bin needs to be modeled, accounting for all case dimensions and orientations. Each case can have two possible orientations considering rotation along a given axis, resulting in a total of six possible orientations for all three axes, as shown in Figure 2.

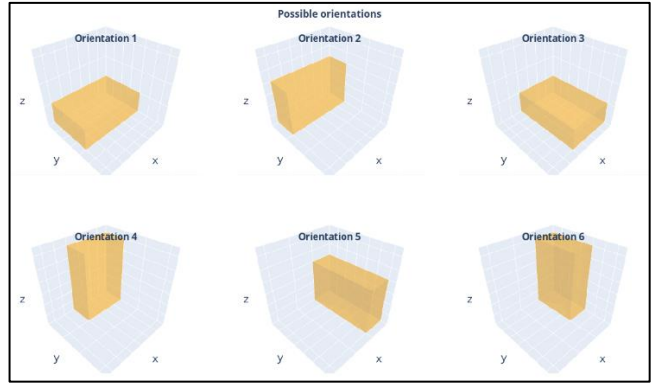

Figure 2. Six possible different orientations of a case.

From a practical standpoint, there can be additional requirements, such as ensuring that each loaded case has sufficient area support.

The MIP formulation presented in the next section is based on a basic model we developed for 3D bin-packing and made available online at https://github.com/dwave-examples/3d-bin-packing. We formally describe that basic model in this paper and further add support constraints to it. We then perform a comprehensive analysis on the test cases we developed using classical solver (Gurobi) and hybrid classical-quantum annealing solver (CQM). It is important to note that the above models were developed without requiring any intermediate steps to translate the problem into a QC friendly form, and that the formulations were directly solvable using the CQM solver as well as commercially available classical solvers such as Gurobi, Xpress, CPLEX, etc., which enormously facilitated the testing and comparison between classical and quantum hybrid solver performance on the same formulation. Appendix III offers an overview of D-Wave's Hybrid Solver Solutions, providing more details on the CQM solver, including the types of variables and constraints, as well as the maximum number of variables and constraints it can handle.

## 2.1 Mathematic Model

In this section, we present the exact MIP formulation for the 3D bin-packing problem. While our results are presented for single bin, we developed our model to handle multi-bin packing problems. An exact model for 3D bin-packing under practical constraints, such as support and fragility constraints, is proposed in Paquay et al. (2016). As part of this study, we propose an alternative formulation that uses six binary variables to model orientation, as opposed to the nine variables used in previous models (Paquay et al., 2016), leading to three fewer binary variables per case for orientation modeling. This reduction in variable count decreases the complexity of the associated geometric, boundary, and support constraints while still capturing all the necessary information.

We use $I$ to represent set of cases, $J$ as set of bins, and $K$ as set of orientation for each case. Each bin ($j$) has three dimensions: width ($W_j$), length ($L_j$), and height ($H_j$). Refer to Appendix I for a full list of sets, parameters, and variables used to model 3D bin-packing problem using the MILP model. To model multiple bins, it is assumed that bins are located back-to-back next to each other along the x-axis of the coordinate system; this way, only one coordinate system with width of $\max_j W_j$, height of $\max_j H_j$, and length of $\sum_j L_j$ to model all bins is required. With this simplification, the first bin is located along the $x$ axis of the bin at $0 \leq x \leq L_1$ interval, the second at $L_1 < x \leq L_1 + L_2$, and the last bin $n$ is located at $\sum_{j=1}^{n-1} L_j < x \leq \sum_{j=1}^{n} L_j$. Constraints are then added to ensure that cases are not placed between any two bins. The following set is used to help simplify the notations for the equations in this section: $\mathcal{L}_j = L_1, L_1 + L_2, \ldots, \sum_{j'=1}^{n} L_{j'}$.

I. Objective

The objective is stated as follows:

$$\text{Min } obj = \frac{\sum_{i \in I}(z_i + z'_i)}{|I|} \frac{l_i\, w_i\, h_i}{\max_{i'}\,(l_{i'}\, w_{i'}\, h_{i'})} + \sum_{j \in J} g_j + \sum_j H_j\, e_j. \tag{1}$$

Where $l_i, w_i,$ and $h_i$ are length, width, and height of case $i$. $z_i$ is the location of the case along $z$ axis of the bin, $z'_i$ is effective height of the case considering rotation of the case (see equation (5)). $g_j$ is the height of the case at the top of the bin $j$, $e_j$ is a binary variable that defines whether bin $j$ is used. $|I|$ defines number of cases.

The first term ensures that each case is packed down. This term has two components: the first component models the average height of the cases in a bin, while the second component normalizes this average by volume. Without the second component, bulky cases are usually packed at the top of the bin and small cases are packed at the bottom.

The second term in the objective ensures that the height of the topmost case in each bin is minimized. This objective is weakly considered in the first objective term, but here is given more importance.

The third term in the objective function minimizes the total number of bins. Note that this term is multiplied by the height of the bins so its contribution to the objective has the same weight as the first and second objective terms.

## II. Orientation Constraints

As shown in Figure 2, each case can be placed in six orientations in a bin. We use the following equation to enforce exactly one orientation per case:

$$\sum_{k \in K} r_{i,k} = 1 \quad \forall i \in I, \tag{2}$$

where $r_{i,k}$ is a binary variable that defines if case $i$ is packed at orientation $k$. Orientation defines the effective length ($x'_i$), width ($y'_i$), and height ($z'_i$) of the cases along $x$, $y$, and $z$ axes.

$$x'_i = l_i\, r_{i,1} + l_i\, r_{i,2} + w_i\, r_{i,3} + w_i\, r_{i,4} + h_i\, r_{i,5} + h_i\, r_{i,6} \quad \forall i \in I \tag{3}$$

$$y'_i = w_i\, r_{i,1} + h_i\, r_{i,2} + l_i\, r_{i,3} + h_i\, r_{i,4} + l_i\, r_{i,5} + w_i\, r_{i,6} \quad \forall i \in I \tag{4}$$

$$z'_i = h_i\, r_{i,1} + w_i\, r_{i,2} + h_i\, r_{i,3} + l_i\, r_{i,4} + w_i\, r_{i,5} + l_i\, r_{i,6} \quad \forall i \in I \tag{5}$$

## III. Case and Cuboid Assignment Constraints

We use $u_{i,j}$ binary variable to define if case $i$ is assigned to bin $j$. We then use the following three constraints to ensure correct assignment of the cases to bins:

$$\sum_{j \in J} u_{i,j} = 1 \quad \forall i \in I \tag{6}$$

$$u_{i,j} \leq e_j \quad \forall j \in J\, \forall i \tag{7}$$

$$e_j \geq e_{j+1} \quad \forall p \in P \mid \forall j \in T_p\ \forall j+1 \in T_p, \tag{8}$$

where $T_p$ is a set containing number of bins of type $p$. Equation (6) ensures that each case can go in to exactly one bin, (7) ensures that cases can only be assigned to bins that are in use. Equation (8) ensures that bins of the same type (same size) are added in order, i.e., if there are multiple numbers of bins of type $p$, bin $j$ of type $p$ is in use before bin $(j + 1)$ of that type.

## IV. Geometric Constraints

Geometric constraints are applied to prevent overlap between cases. In the following constraints, "left" and "right" refer to the position of the case along the $x$ axis of a bin, "behind" and "in front of" to the $y$ axis, and "above" and "below" to the $z$ axis. To avoid overlap between each pair of cases, at least one of the following must occur:

- case $i$ is on the left of case $i'$ ($q = 0$):

$$-\left(2 - u_{i,j}\, u_{i',j} - b_{i,i',0}\right)\mathcal{L}_n + x_i + x_i' - x_{i'} \leq 0 \quad \forall i, i' \in I \;\; \forall j \in J \tag{9}$$

- case $i$ is behind case $i'$ ($q = 1$):

$$-\left(2 - u_{i,j}\, u_{i',j} - b_{i,i',1}\right)\mathcal{W} + y_i + y_i' - y_{i'} \leq 0 \quad \forall i, i' \in I \;\forall j \in J \tag{10}$$

- case $i$ is below case $i'$ ($q = 2$):

$$-\left(2 - u_{i,j}\, u_{i',j} - b_{i,i',2}\right)\mathcal{H} + z_i + z_i' - z_{i'} \leq 0 \quad \forall i, i' \in I \;\forall j \in J \tag{11}$$

- case $i$ is on the right of case $i'$ ($q = 3$):

$$-\left(2 - u_{i,j}\, u_{i',j} - b_{i,i',3}\right)\mathcal{L}_n + x_{i'} + x'_{i'} - x_i \leq 0 \quad \forall i, i' \in I \;\forall j \in J \tag{12}$$

- case $i$ is in front of case $i'$ ($q = 4$):

$$-\left(2 - u_{i,j}\, u_{i',j} - b_{i,i',4}\right)\mathcal{W} + y_{i'i} + y'_{i'} - y_i \leq 0 \quad \forall i, i' \in I \;\forall j \in J \tag{13}$$

- case $i$ is above case $i'$ ($q = 5$):

$$-\left(2 - u_{i,j}\, u_{i',j} - b_{i,i',5}\right)\mathcal{H} + z_{i'} + z'_{i'} - z_i \leq 0 \quad \forall i, i' \in I \;\forall j \in J. \tag{14}$$

The following is then required to enforce only one of the above constraints:

$$\sum_{q \in Q} b_{i,i',q} = 1 \quad \forall i', i \in I, \tag{15}$$

where binary variable $b_{i,i',q}$ is used to define geometric relation $q$ between cases $i$ and $i'$.

## V. Cuboid Boundary Constraints

The sets of constraints below ensures that case $i$ is not placed outside the bins. When $u_{i,j}$ is 0, these constraints are relaxed.

$$x_i + x'_i - \mathcal{L}_j \leq (1 - u_{i,j}) \mathcal{L}_n \quad \forall j \in J \; \forall i \in I \tag{16}$$

$$x_i - u_{i,j} \mathcal{L}_j \geq 0 \quad \forall j \in J \quad \forall i \in I \tag{17}$$

$$y_i + y'_i - W_j \leq (1 - u_{i,j}) \mathcal{W} \quad \forall j \in J \; \forall i \in I \tag{18}$$

$$z_i + z'_i - H_j \leq (1 - u_{i,j}) \mathcal{H} \quad \forall j \in J \; \forall i \in I \tag{19}$$

$$z_i + z'_i - g_j \leq (1 - u_{i,j}) \mathcal{H} \quad \forall j \in J \; \forall i \in I, \tag{20}$$

where $\mathcal{W}, \mathcal{H}$ are max width and height across all bins. $\mathcal{L}_j$ defines end location of each bin when they are set back-to-back next to each other along the $x$ axis and defined as $\mathcal{L}_j = \sum_{j'=1}^{j} L_{j'}$.

## VI. Support Constraints

For any two boxes $i$ and $i'$, $s_{i,i'}$, is used as a continuous variable to define amount of support that case $i'$ provides to case $i$. The binary variable $f_{i,i'}$ is used to determine if case $i$ and $i'$ touch each other along the $z$ axis.

The following equation ensures that each case receives required support from all other cases.

$$\sum_{i' \in I \;|\; i \neq i'} s_{i,i'} \geq T x'_i y'_i \quad \forall \, i \in I \tag{21}$$

Where $0 \leq T \leq 1$ is support threshold, $x'_i$ is the effective length of case $i$, and $y'_i$ is its effective width.

Considering equation (2), any products of $r_{i,\cdot}$ variables are zero, which results in this linear expression for $x'_i\, y'_i$ term.

$$x'_i\, y'_i = l_i\, w_i\, (r_{i,1} + r_{i,3}) + l_i\, h_i(r_{i,2} + r_{i,5}) + w_i\, h_i(r_{i,4} + r_{i,6})$$

The binary variable $f_{i,i'}$ is used to determine if case $i$ and $i'$ touch each other along the $z$ axis of the bin and if case $i'$ provides support for case $i$.

$$-\mathcal{H}\,(1 - f_{i,i'}) + z_{i'} + z'_{i'} \leq z_i \leq z_{i'} + z'_{i'} + \mathcal{H}\,(1 - f_{i,i'}) \quad \forall i, i' \tag{22}$$

$$x_{i'} \leq x_i + x'_i + (1 - f_{i,i'})\mathcal{L}_n \quad \forall i, i' \tag{23}$$

$$(1 - f_{i,i'})\mathcal{L}_n + x_{i'} + x'_{i'} \geq x_i \quad \forall i, i' \tag{24}$$

$$y_{i'} \leq y_i + y'_i + \mathcal{W}\,(1 - f_{i,i'}) \quad \forall i, i' \tag{25}$$

$$(1 - b_{i,i'})\mathcal{W} + y_{i'} + y'_i \geq y_i \quad \forall i, i' \tag{26}$$

When $f_{i,i'} = 1$, the above constraints are enforced, when $f_{i,i'} = 0$, they are relaxed.

If the two cases $i$ and $i'$ touch each other along $z$ axis of the bin, $f_{i,i'} = 1$, the amount of support that case $i'$ provides to case $i$, $s_{i,i'}$, is nonzero and defined as product of overlaps along $x$ and $y$ axis ($o^x_{i,i'}$ and $o^y_{i,i'}$,) of the bin.

$$s_{i,i'} = f_{i,i'}\, o^x_{i,i'}\, o^y_{i,i'}, \tag{27}$$

In Appendix III, we demonstrate the definition of these overlapping variables and present a method to linearize the nonlinear equation in (27).

## 3. Results Summary

In this section, we present key results from the computational comparison between the Gurobi and CQM solvers on the exact formulations presented earlier. To test and compare purely classical and quantum hybrid methods for the 3D bin-packing problem, extensive computational testing was conducted on a testbed of instances based on sanitized data from real-world business use cases. Appendix IV provides information for each use case including sizes and quantities of all cases and bins.

The solvers used in the study were Gurobi 11.0.3. for classical computation and D-Wave's CQM solver 1.12 for quantum hybrid solution. Gurobi solver was run on an Intel® Xeon® Platinum 8272CL CPU @ 2.60GHz using 8 threads. The hardware for the CQM solver was not provided by D-Wave. Gurobi offers a wide range of solver parameters, whereas CQM does not provide any. While having numerous solver parameters adds flexibility in problem-solving, adjusting those for optimal solver performance is often a complex task. We used Gurobi with its default solver parameters.

To evaluate solution quality, we conducted tests with progressively increasing time limits, recording the relative gap and bin utilization at the conclusion of each time limit for every problem instance. The relative gap measures the difference between the best known integer feasible solution relative to its *Best Bound* solution.

$$Relative\ Gap\ \% = \frac{Best\ Integer\ Solution\ -\ Best\ Bound}{Best\ Integer\ Solution}$$

The Best Bound solution is the value of the objective function where some integer variables are fixed at their best possible integer values, while others (those that have not yet been branched on or fully explored) are relaxed to continuous variables. We use the same Best Bound, provided by Gurobi, to evaluate the *Relative Gap* for both solvers.

Bin utilization is the total volume of cases placed in the bin divided by the volume of the bin up to the maximum occupied height in the solution. It naturally implies that higher utilization results in superior solution quality.

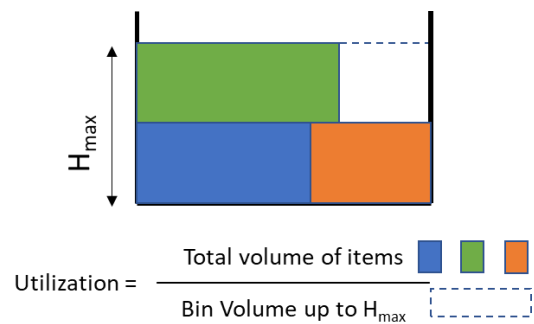

Figure 4. Bin utilization defined as the percentage obtained by dividing the volume of packed cases by the available volume of the occupied bin.

In the following section, we present computational results comparing the Gurobi and CQM solvers on the MIP model with and without support constraints. To model specific use cases and ensure tractability, the height of the cases is assumed fixed and six orientations along the vertical axis are allowed. We tested instances containing up to 158 cases for a preset time limit of 25, 50, 100, 200, 400, and 800 seconds. Relative gap and utilization values are calculated when there is at least one feasible solution.

### 3.1. Results Without Support Constraints

In the first set of experiments on the MIP model, we did not include support area constraints. Figure 5 illustrates the relative gap over time for each instance. Both solvers generate a solution pool (referred to as a sample set for the CQM solver). For the CQM solver, the solutions are distinct, enabling statistical analysis and the creation of graphs for the minimum, median, $25^{th}$, and 75th percentiles. In contrast, although Gurobi provided multiple solutions, they were identical, preventing the generation of statistical curves.

The results indicate that Gurobi tends to perform better than the CQM solver for smaller problem sizes when enough time is provided. As problem size increases, however, Gurobi loses its advantage and is not able to produce any feasible solution for some of the large test cases. The CQM solver produces feasible solutions across all problem sets and performs better than Gurobi when the problem size is larger or when the time limit is short for smaller problem sizes.

Utilizations metrics carry the same information as relative gap about performance of the solvers. For interested readers, we have included utilization plots in Appendix V.

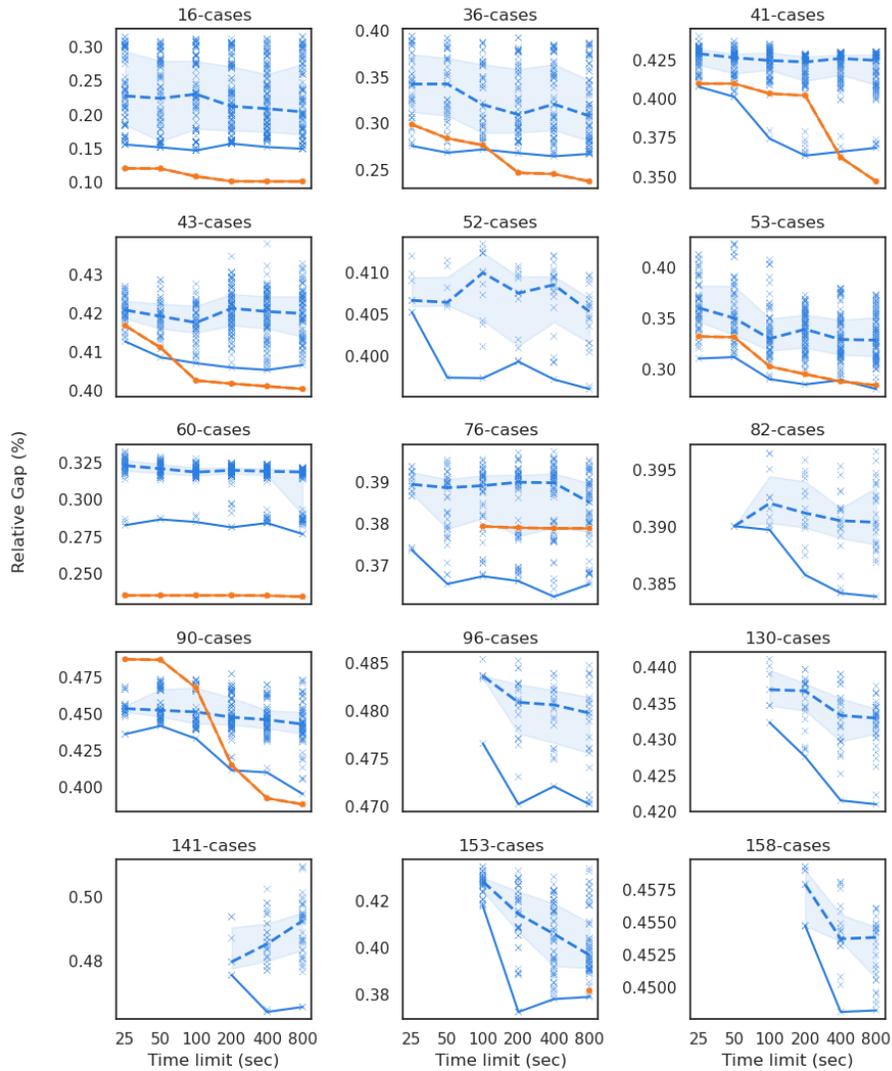

Figure 5: Comparison of the Gurobi and CQM solvers on the MIP model (without support constraints) based on relative gap (%) versus time limit (sec) across all 3D bin-packing problem instances. Blue curves represent the CQM solver, while orange curves represent Gurobi. The dashed blue line indicates the median values for the CQM solver's solutions, and the solid blue line shows the minimum value among its solutions. Blue "x" marks indicate individual CQM solver solutions. The shaded blue area represents the 25th and 75th percentiles around the median curve. The median, 25th, and 75th percentiles are not provided for Gurobi as its solution is deterministic. Missing curves or data points indicate that the solver did not provide a feasible solution within the given time limit.

## 3.2. Results with Support Constraints

In our final test, we introduced constraints for area support (see section 2.1.VI) in combination with constraints for the base model. Support constraints not only extend the problem size but also introduce significant complexity to the model. In this test, we set the support threshold at 80%; that is, we set $T$ to 0.8 in (21). Gurobi managed to return a feasible solution for only the smallest problem instance (16-cases), while the CQM solver returned feasible solutions for two instances (16-and 41-cases), with both solvers failing to provide feasible solutions for the remaining test cases. Similar to previous results, Gurobi achieved lower relative gap and higher utilization for the small test case, while the CQM solver was able to produce feasible solutions for more problem instances. The results for relative gap are graphically shown in Figure 6 and plot for utilization is included in Appendix V.

This test reaffirms the observation that the Gurobi solver, while efficient at handling smaller problems, tends to lose its advantage against the CQM solver as the problem size increases and/or more complex problem features are taken into consideration.

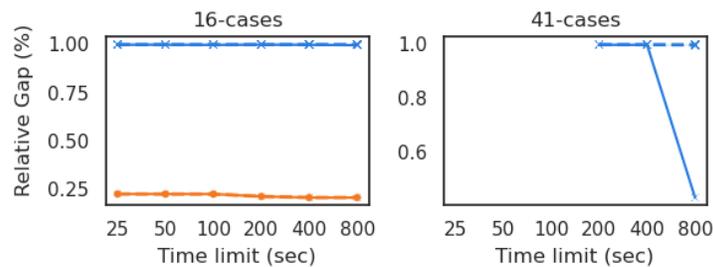

Figure 6: fig (a) relative gap vs time limit (sec) for the CQM vs Gurobi solvers with support constraints. Blue curves represent the CQM solver, while the orange curves represent Gurobi.

## 4. Key Takeaways and Future Outlook

Optimization problems entail a common fundamental principle: finding the best solution to maximize or minimize one or more objectives while satisfying certain constraints. They are hard to solve since the space of feasible solutions expands very quickly, often exponentially, making it increasingly challenging to traverse through the solution space in an efficient manner. In view of this, it is desirable to investigate and develop solution methods whose performances scale well with an increase in problem size and/or as more complex problem features are taken into consideration. In this paper, we have explored the use of annealing-based hybrid quantum computing, framed as the CQM solver, to solve a business optimization problem in comparison with a commercially available classical optimization solver that is traditionally used to solve such problems.

From the computational testing on the 3D bin-packing problem, it is concluded that the CQM solver can solve conventionally hard optimization problems, and in many instances, produce results that are competitive with the best commercially available classical solvers. A key distinguishing feature of the CQM solver architecture is that the original problem formulation (for example, an MILP or MIQP) can be directly fed into the CQM solver without any requirements to translate it into a QC-friendly form.

For the base model without support constraints, Gurobi performed better on small instances, especially when the model is allowed to run for a sufficiently long time. As the problem size increased, however, the CQM solver was found to surpass Gurobi in the likelihood of finding a feasible solution.

A practical 3D bin-packing problem requires support constraints to develop usable solutions, avoiding instances of floating cases with insufficient support. In this study, both the Gurobi and CQM solvers struggled to generate good feasible solutions for most problem instances with support constraints. The classical solver managed to solve only one instance, while the hybrid classical-quantum solver successfully solved two.

Researchers worldwide are trying to build smarter quantum computers that can surpass the computing power of the most powerful classical supercomputers and routinely achieve quantum advantage in real-world optimization problems. As the next generation of quantum hardware continues to advance and quantum hybrid solution approaches evolve, hybrid solvers have potential to become potent tools for tackling a wider range of combinatorial optimization problems, bridging the gap with classical solvers, and potentially surpassing them in certain domains, which showcases the exciting potential of quantum computing in revolutionizing the field of optimization.

# Appendix I: D-Wave's Hybrid Solver Solutions

The cloud-based hybrid solver service (HSS) leverages classical and quantum solution methods to solve optimization problems that are too large for QPUs using a portfolio of heuristic solvers.

In the hybrid workflow, classical and quantum solution methods work in tandem that leads to convergence to good solutions faster than purely classical workflow, a phenomenon referred to as *hybrid acceleration*. The HSS supports two different types of problem formulations (or models) to solve unconstrained and constrained quadratic optimization problems: binary and constrained quadratic models. The QUBO formulation used by the binary quadratic model (BQM) solver is the native formulation of all D-Wave QPUs and, consequently, this solver shows the best solution performance in solving QUBOs. In principle, any NP-hard problem can be translated to a QUBO formulation, though in actual practice, this poses a challenge for some problems and the solver performance is impacted by subpar translation quality. To encompass optimization models that do not match the native formulation of the quantum hardware, the CQM solver provides interface support to end users for a much broader set of business optimization applications. Unlike the BQM solver, the CQM solver supports the optimization constraints to be incorporated as hard constraints, thus allowing for more realistic representation of real-world problems by avoiding strictly infeasible regions of the solution space. One limitation of the CQM solver is that it currently does not support the representation of optimization models on continuous decision variables for certain types of constraints, such as handling quadratic interactions (product) between continuous (non-integer real-valued) decision variables.

Based on solver type and hardware architecture, there are specific limits on the problem size that can be solved using each solver. These limits are specified in Table 1 (D-Wave, 2022).

|  | BQM | CQM |
|---|---|---|
| Objective Function | linear & quadratic | linear & quadratic |
| Variable Type | binary | binary, integer, real |
| Max Values per Variable | 2 | 2, +/- $2^{53}$ |
| Constraint Representation | via penalties | variable bounds Real/integer linear & quadratic equality/inequality |
| Max Variables | 1 million | 500,000 |
| Max Constraints | – | 100,000 |
| Max Biases | 200 million | 2 billion |

Table A. 1. Hardware limits for BQM and CQM solvers in D-Wave's hybrid solver service (as of Jan 2024).

## Appendix II: Nomenclature

Below is a list of constants, variables and parameters that are used in the MIP model.

### AII.1 Parameters and Sets

| |
|---|
| $P$: set of bins with p as index defining type of the bin |
| $T_p$: set containing number of bins of type $p$ |
| $J$: ordered set of bins |
| $n$: number of bins |
| $I$: set of cases |
| $m$: number of cases |
| $K$: possible geometric relation (e.g., behind, above, etc.) between every pair of cases  0,1, … ,5 |
| $L_j$: length of bin |
| $\mathcal{L} = L_1, L_1 + L_2, \dots, \sum_{j'=1}^{n} L_{j'}$ : end location of each bin when they are set back-to-back next to each other along the $x$ axis. |
| $W_j$: width of bin $j$ |
| $\mathcal{W}$: max width across all bins |
| $H_j$: height of the bin $j$ |
| $\mathcal{H}$: max height across all bins |
| $l_i$: height of case $i$ |
| $w_i$: width of case $i$ |
| $h_i$: Height of the case $i$ |
| $T$: support threshold in range [0, 1] |
| $A_{i,i'}$: is the maximum overlap area between cases i and i'. |

## AII.2 Variables

| |
|---|
| $e_j$: binary variable that shows if bin $j$ is used |
| $u_{i,j}$: binary variable that shows if case $i$ is added to bin $j$ |
| $b_{i,i',q}$: binary variable defining geometric relation $q$ between cases $i$ and $i'$ |
| $g_j$: continuous variable showing height of the topmost case in bin $j$ |
| $r_{i,k}$: are binary variables defining $k$ orientations for case $i$ |
| $x_i, y_i, z_i$: continuous variable defining location of the back lower left corner of case $i$ along $x$, $y$, and $z$ axes of the bin |
| $x'_i, y'_i, z'_i$: effective length, width, and height of box $i$, considering orientation, along $x$, $y$, and $z$ axis of the bin |
| $s_{i,i'}$: continuous variable, to define amount of support that box $i'$ provides to box $i$ |
| $f_{i,i'}$: determines if box $i$ and $i'$ touch each other along $z$ axis |
| $s_{-1,i}$: continuous variable, to define amount of support that ground provides to box $i$. |
| $f_{i,-1}$: determines if box $i$ is placed at the bottom of the bin |
| $o^x_{i,i'}$: overlap between case $i$ and $i'$ along $x$ axis |
| $o^y_{i,i'}$: overlap between case $i$ and $i'$ along $y$ axis |
| $o^{xy}_{i,i'}$: overlap between case $i$ and $i'$ in $xy$ surface |

## Appendix III: Linearization of Support Constraints

The support equation in (27) ($s_{i,i'} = f_{i,i'} \ o^x_{i,i'} \ o^y_{i,i'}$) contains multiplications of three variables. This constraint can be turned into linear and quadratic constraints using the following inequalities:

$$s_{i,i'} \leq A_{i,i'} \ f_{i,i'} \tag{A1}$$

$$s_{i,i'} \leq o^x_{i,\,i'} o^y_{i,\,i'}. \tag{A2}$$

Where $A_{i,i'}$ is the maximum possible area of case $i$ and $i'$ in any possible orientations.

Equation A2 involves quadratic interactions among continuous variables. Certain optimization solvers can automatically manage these quadratic interactions. Alternatively, the Piecewise McCormick transformation, as outlined in Castro (2015), can be employed to linearize this equation if necessary. Another strategy, particularly suitable for CQM solver, involves converting continuous variables into integer variables with a specified precision level and using the product of these integer variables to represent the equation.

$o^x_{i,i'}$ and $o^y_{i,i'}$, overlaps along $x$ and $y$ axis, are defined as follows:

$$o^x_{i,\,i'} = \min(x_i + x'_i - x_{i'}, x_{i'} + x'_{i'} - x_i, x'_{i'}, x'_i)$$

$$o^y_{i,\,i'} = \min(y_i + y'_i - y_{i'}, y_{i'} + y'_{i'} - y_i, y'_i, y'_{i'})$$

The minimum function above is defined using these sets of linear equations.

$$o^x_{i,\,i'} \leq x_i + x'_i - x_{i'} + (1 - f_{i,\,i'}) \mathcal{L}_n \tag{A3}$$

$$o^x_{i,\,i'} \leq x_{i'} + x'_{i'} - x_i + (1 - f_{i,\,i'}) \mathcal{L}_n \tag{A4}$$

$$o^x_{i,i'} \leq x'_i \tag{A5}$$

$$o^x_{i,i'} \leq x'_{i'} \tag{A6}$$

$$o^y_{i,i'} \leq y_i + y'_i - y_{i'} + (1 - f_{i,i'})\mathcal{W} \tag{A7}$$

$$o^y_{i,i'} \leq y_{i'} + y'_{i'} - y_i + (1 - f_{i,i'})\mathcal{W} \tag{A8}$$

$$o^y_{i,i'} \leq y'_i \tag{A9}$$

$$o^y_{i,i'} \leq y'_{i'} \tag{A10}$$

Finally, cases can also get support from the ground if they are placed at the bottom of the bin.

We use $s_{-1,i}$ to indicate the amount of support that the ground provides to a given case.

$$z_i \leq (1 - f_{-1,i})\mathcal{H} \quad \forall i \tag{A11}$$

$$s_{-1,i} \leq f_{-1,i} A_i \quad \forall i \tag{A12}$$

Where $A_i$ is the maximum area of case $i$ in any possible orientation.

## Appendix IV: Use Case Information

| Instance ID | Number of cases | Case Description | | | | | Bin Description | | | |
|---|---|---|---|---|---|---|---|---|---|---|
| | | IDs | Quantity | Length | Width | Height | Quantity | Length | Width | Height |
| 1 | 16 | 0 | 1 | 10.88 | 9.82 | 10.87 | 1 | 50.00 | 50.00 | 50.00 |
| | | 1 | 1 | 10.88 | 9.82 | 14.13 | | | | |
| | | 2 | 1 | 10.88 | 15.98 | 11.49 | | | | |
| | | 3 | 1 | 10.88 | 15.98 | 13.51 | | | | |
| | | 4 | 1 | 11.85 | 24.2 | 9.43 | | | | |
| | | 5 | 1 | 11.85 | 24.2 | 15.57 | | | | |
| | | 6 | 1 | 14.55 | 24.2 | 10.86 | | | | |
| | | 7 | 1 | 14.55 | 24.2 | 14.14 | | | | |
| | | 8 | 1 | 15.52 | 12.39 | 9.88 | | | | |
| | | 9 | 1 | 15.52 | 12.39 | 15.12 | | | | |
| | | 10 | 1 | 15.52 | 13.41 | 11.44 | | | | |
| | | 11 | 1 | 15.52 | 13.41 | 13.56 | | | | |
| | | 12 | 1 | 23.6 | 12.61 | 11.74 | | | | |
| | | 13 | 1 | 23.6 | 12.61 | 13.26 | | | | |
| | | 14 | 1 | 23.6 | 16.53 | 25 | | | | |
| | | 15 | 1 | 23.6 | 20.86 | 25 | | | | |
| 2 | 36 | IDs | Quantity | Length | Width | Height | Quantity | Length | Width | Height |
| | | 0 | 1 | 8.90 | 7.45 | 11.09 | 1 | 50.00 | 50.00 | 50.00 |
| | | 1 | 1 | 8.90 | 9.08 | 11.09 | | | | |
| | | 2 | 1 | 9.79 | 12.61 | 11.74 | | | | |
| | | 3 | 1 | 9.94 | 12.61 | 13.26 | | | | |

| | | | | | | | | | |
|---|---|---|---|---|---|---|---|---|---|
| | | 4 | 1 | 9.95 | 7.97 | 12.29 | | | |
| | | 5 | 1 | 9.95 | 12.89 | 12.29 | | | |
| | | 6 | 1 | 10.88 | 9.82 | 10.87 | | | |
| | | 7 | 1 | 10.88 | 9.82 | 14.13 | | | |
| | | 8 | 1 | 10.88 | 15.98 | 11.49 | | | |
| | | 9 | 1 | 10.88 | 15.98 | 13.51 | | | |
| | | 10 | 1 | 11.34 | 10.19 | 12.71 | | | |
| | | 11 | 1 | 11.34 | 10.67 | 12.71 | | | |
| | | 12 | 1 | 11.34 | 8.00 | 13.91 | | | |
| | | 13 | 1 | 11.34 | 8.53 | 13.91 | | | |
| | | 14 | 1 | 11.85 | 9.52 | 15.57 | | | |
| | | 15 | 1 | 11.85 | 9.79 | 9.43 | | | |
| | | 16 | 1 | 11.85 | 14.41 | 9.43 | | | |
| | | 17 | 1 | 11.85 | 14.68 | 15.57 | | | |
| | | 18 | 1 | 12.26 | 7.33 | 13.91 | | | |
| | | 19 | 1 | 12.26 | 9.20 | 13.91 | | | |
| | | 20 | 1 | 12.26 | 8.37 | 12.71 | | | |
| | | 21 | 1 | 12.26 | 12.49 | 12.71 | | | |
| | | 22 | 1 | 13.65 | 8.78 | 12.29 | | | |
| | | 23 | 1 | 13.65 | 12.08 | 12.29 | | | |
| | | 24 | 1 | 13.66 | 12.61 | 13.26 | | | |
| | | 25 | 1 | 13.81 | 12.61 | 11.74 | | | |
| | | 26 | 1 | 14.55 | 9.68 | 10.86 | | | |
| | | 27 | 1 | 14.55 | 10.32 | 14.14 | | | |
| | | 28 | 1 | 14.55 | 13.88 | 14.14 | | | |
| | | 29 | 1 | 14.55 | 14.52 | 10.86 | | | |
| | | 30 | 1 | 14.70 | 6.21 | 11.09 | | | |
| | | 31 | 1 | 14.70 | 10.32 | 11.09 | | | |
| | | 32 | 1 | 15.52 | 12.39 | 9.88 | | | |
| | | 33 | 1 | 15.52 | 12.39 | 15.12 | | | |
| | | 34 | 1 | 15.52 | 13.41 | 11.44 | | | |
| | | 35 | 1 | 15.52 | 13.41 | 13.56 | | | |
| 3 | 41 | IDs | Quantity | Length | Width | Height | Quantity | Length | Width | Height |
| | | 0 | 32 | 4.30 | 8.00 | 8.10 | 1 | 38.10 | 38.10 | 22.00 |
| | | 1 | 9 | 20.00 | 4.00 | 13.50 | | | | |
| 4 | 43 | IDs | Quantity | Length | Width | Height | Quantity | Length | Width | Height |
| | | 0 | 12 | 14.60 | 5.70 | 5.80 | 1 | 38.10 | 38.10 | 22.00 |
| | | 1 | 21 | 11.00 | 5.00 | 10.00 | | | | |
| | | 2 | 10 | 16.10 | 6.60 | 5.00 | | | | |
| 5 | 52 | IDs | Quantity | Length | Width | Height | Quantity | Length | Width | Height |
| | | 0 | 32 | 4.60 | 9.20 | 8.70 | 1 | 38.10 | 38.10 | 22.00 |
| | | 1 | 20 | 16.20 | 6.20 | 6.00 | | | | |
| 6 | 53 | IDs | Quantity | Length | Width | Height | Quantity | Length | Width | Height |
| | | 0 | 3 | 41.40 | 24.40 | 24.80 | 1 | 116.50 | 116.50 | 130.00 |
| | | 1 | 2 | 41.50 | 23.00 | 14.00 | | | | |
| | | 2 | 3 | 23.80 | 25.00 | 19.80 | | | | |
| | | 3 | 2 | 28.00 | 14.50 | 22.00 | | | | |
| | | 4 | 8 | 26.00 | 38.20 | 16.50 | | | | |
| | | 5 | 8 | 38.20 | 26.00 | 16.50 | | | | |
| | | 6 | 2 | 37.50 | 23.50 | 29.50 | | | | |
| | | 7 | 1 | 44.80 | 16.20 | 24.30 | | | | |
| | | 8 | 2 | 28.50 | 28.90 | 22.90 | | | | |
| | | 9 | 6 | 32.80 | 24.20 | 17.70 | | | | |
| | | 10 | 6 | 12.20 | 13.50 | 20.00 | | | | |
| | | 11 | 1 | 34.00 | 28.50 | 40.50 | | | | |
| | | 12 | 3 | 38.50 | 26.50 | 16.80 | | | | |
| | | 13 | 2 | 45.50 | 37.00 | 16.00 | | | | |
| | | 14 | 4 | 37.00 | 29.50 | 18.50 | | | | |
| 7 | 60 | IDs | Quantity | Length | Width | Height | Quantity | Length | Width | Height |
| | | 0 | 40 | 8.00 | 5.00 | 5.00 | 1 | 40.00 | 40.00 | 20.00 |
| | | 1 | 20 | 10.00 | 8.00 | 8.00 | | | | |
| 8 | 76 | IDs | Quantity | Length | Width | Height | Quantity | Length | Width | Height |
| | | 0 | 40 | 7.00 | 7.00 | 6.10 | 1 | 38.10 | 38.10 | 25.00 |

|   |   |   | 1 | 36 | 8.00 | 7.00 | 7.20 |   |   |   |   |
|---|---|---|---|---|---|---|---|---|---|---|---|
| 9 | 82 | IDs | Quantity | Length | Width | Height | Quantity | Length | Width | Height |
|   |   | 0 | 30 | 7.60 | 7.60 | 8.50 | 1 | 38.10 | 38.10 | 25.00 |
|   |   | 1 | 28 | 5.80 | 5.10 | 10.00 |   |   |   |   |
|   |   | 2 | 24 | 15.80 | 4.50 | 3.20 |   |   |   |   |
| 10 | 90 | IDs | Quantity | Length | Width | Height | Quantity | Length | Width | Height |
|   |   | 0 | 16 | 14.60 | 5.70 | 5.80 | 1 | 38.10 | 38.10 | 40.00 |
|   |   | 1 | 56 | 5.50 | 5.00 | 10.00 |   |   |   |   |
|   |   | 2 | 18 | 16.10 | 5.50 | 5.00 |   |   |   |   |
| 11 | 96 | IDs | Quantity | Length | Width | Height | Quantity | Length | Width | Height |
|   |   | 0 | 36 | 3.50 | 3.50 | 6.10 | 1 | 21.91 | 22.50 | 25.00 |
|   |   | 1 | 30 | 4.00 | 3.50 | 7.20 |   |   |   |   |
|   |   | 2 | 30 | 3.70 | 4.20 | 7.70 |   |   |   |   |
| 12 | 130 | IDs | Quantity | Length | Width | Height | Quantity | Length | Width | Height |
|   |   | 0 | 48 | 7.60 | 7.60 | 8.50 | 1 | 48.20 | 38.90 | 35.40 |
|   |   | 1 | 32 | 5.80 | 5.10 | 10.00 |   |   |   |   |
|   |   | 2 | 50 | 15.80 | 4.50 | 3.20 |   |   |   |   |
| 13 | 141 | IDs | Quantity | Length | Width | Height | Quantity | Length | Width | Height |
|   |   | 0 | 33 | 14.60 | 5.70 | 6.00 | 1 | 48.20 | 38.90 | 60.00 |
|   |   | 1 | 51 | 11.00 | 5.00 | 6.00 |   |   |   |   |
|   |   | 2 | 32 | 16.10 | 5.50 | 5.00 |   |   |   |   |
|   |   | 3 | 25 | 16.10 | 6.60 | 5.00 |   |   |   |   |
| 14 | 153 | IDs | Quantity | Length | Width | Height | Quantity | Length | Width | Height |
|   |   | 0 | 57 | 7 | 7 | 6.1 | 1 | 48.2 | 38.9 | 35.4 |
|   |   | 1 | 66 | 4 | 3.5 | 7.2 |   |   |   |   |
|   |   | 2 | 30 | 8 | 7 | 7.2 |   |   |   |   |
| 15 | 158 | IDs | Quantity | Length | Width | Height | Quantity | Length | Width | Height |
|   |   | 0 | 100 | 3.8 | 3.8 | 8.5 | 1 | 38.1 | 38.1 | 30 |
|   |   | 1 | 42 | 5.8 | 5.1 | 10 |   |   |   |   |
|   |   | 2 | 16 | 15.8 | 4.5 | 3.2 |   |   |   |   |

Table A. 2. Description of instances used in the computational study.

## Appendix V: Bin Utilization

Figure A1 shows utilization for the best possible solution from each solver for every 3D bin-packing use case. Utilization is not directly captured in the objective function; however, it is an important metric to evaluate the quality of the 3Dbin-packing solutions.

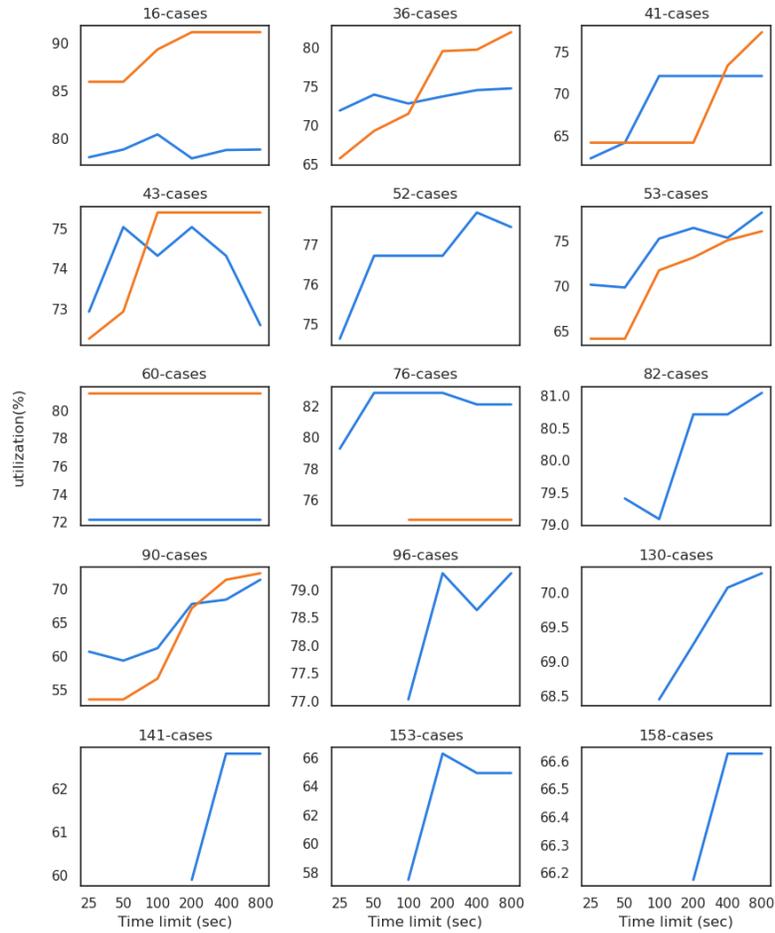

Figure A1: Utilization (%) vs time limit (sec) for 3D-bin packing use cases using Gurobi and CQM solvers. Blue curves represent CQM, while orange curves represent Gurobi. Missing curves indicate that the solver did not provide a feasible solution within the given time limit.

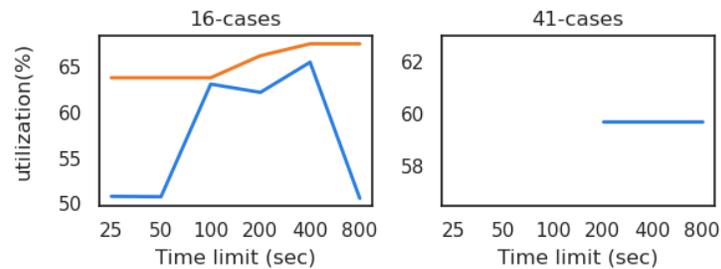

Figure A2: Fig(b) utilization vs time limit (sec) for CQM vs Gurobi with support constraints. Blue curves represent CQM, while the orange curve represents Gurobi.